\documentclass[12pt]{article}
\usepackage{amssymb}

\usepackage{amsmath,latexsym,amsfonts}
\usepackage{graphicx}
\usepackage{srcltx}
\usepackage{amscd}

\newtheorem{lemma}{Lemma}[section]
\newtheorem{coro}[lemma]{Corollary}
\newtheorem{prop}[lemma]{Proposition}
\newtheorem{thm}[lemma]{Theorem}
\newtheorem{defn}[lemma]{Definition}
\makeatletter\@addtoreset{equation}{section}
\renewcommand\theequation{\thesection.\@arabic\c@equation}

\begin{document}
\begin{center}
{\LARGE   Transverse conformal Killing forms on K\"ahler foliations}

 \renewcommand{\thefootnote}{}
\footnote{2000 \textit {Mathematics Subject Classification.}
53C12, 53C27, 57R30}\footnote{\textit{Key words and phrases.}
Transverse Killing form, Transverse conformal Killing form}
\renewcommand{\thefootnote}{\arabic{footnote}}
\setcounter{footnote}{0}

\vspace{1 cm} {\large Seoung Dal Jung}
\end{center}
\vspace{0.5cm}

{\bf Abstract.} On a closed, connected Riemannian manifold with a  K\"ahler foliation of codimension $q=2m$, any transverse Killing $r$ $(\geq 2)$-form is parallel (Jung and Jung, 2012). In this paper, we study transverse conformal Killing forms on K\"ahler foliations. In fact, if the foliation is minimal, then for any transverse conformal Killing $r$-form $\phi$ $(r\ne m, \ 2\leq r\leq q-2)$, $J\phi$ is parallel. Here $J$ is defined in section 4.

\section{Introduction}
On Riemannian manifolds, conformal Killing forms are generalizations of conformal Killing fields, which were introduced by K. Yano [\ref{Yano}] and T.
Kashiwada [\ref{Kashiwada},\ref{Kashiwada1}].
Many researchers have studied the conformal Killing forms [\ref{MS}, \ref{SE}, \ref{TA}, \ref{TA1}].   
On a foliated Riemannian manifold, we can study the analoguous problems. 
Let
$\mathcal F$ be a transversally oriented Riemannian foliation on a
compact oriented Riemannian manifold $M$ with codimension $q$. 
  A
transversal conformal Killing field is a normal field with a flow
preserving the conformal class of the transverse metric. As a generalization of a transversal conformal Killing field, we define  the {\it transverse conformal Killing $r$-forms} $\phi$ as follows: for any
vector field $X$ normal to the foliation,
\begin{align*}
\nabla_X\phi-{1\over r+1}i(X)d\phi+{1\over
q-r+1}X^{\frak b}\wedge\delta_T\phi=0,
\end{align*}
where $r$ is the degree of the form $\phi$ and $X^{\frak b}$ is the dual 1-form
of $X$. For the definition of $\delta_T$,  see Section 3. The transverse conformal Killing
forms $\phi$ with $\delta_T\phi=0$ are called {\it transverse Killing forms}. Recently, S. D. Jung and K. Richardson [\ref{JK}] studied the transverse Killing and conformal Killing forms on Riemannian foliations. And S. D. Jung and M. J. Jung [\ref{JJ2}] studied some properties of the transverse Killing forms on K\"ahler foliations. That is, on a closed, connected Riemannian manifold with a K\"ahler foliation of codimension $q=2m$, any transverse Killing $r (\geq 2)$-form is parallel. 
In this paper, we study the transverse conformal Killing forms on K\"ahler foliations. In section 2, we review  the basic facts on a Riemannian foliation.  In section 3, we study the transverse conformal Killing forms  and curvature properties on Riemannian foliations. In section 4, we study the curvatures  and several operators on  K\"ahler foliations. In section 5, we prove the following: on a K\"ahler foliation with $q=2m$, if $\phi$ is a transverse conformal Killing $m$-form, then $J\phi$ is parallel. In particular, when $(\mathcal F,J)$ is minimal, for any transverse conformal Killing $r$ $(2\leq r\leq q-2)$-forms $\phi$, $J\phi$ is also parallel. Here $J$ is an extension of the complex structure $J$ to the basic forms.

 \section{Preliminaries}
 Let $(M,g_M,\mathcal F)$ be a $(p+q)$-dimensional
Riemannian manifold with a foliation $\mathcal F$ of codimension
$q$ and a bundle-like metric $g_M$ with respect to $\mathcal F$.
 Then there exists an exact sequence of vector bundles
\begin{align}\label{eq1-1}
 0 \longrightarrow L \longrightarrow
TM {\overset\pi\longrightarrow} Q \longrightarrow 0,
\end{align}
where $L$ is the tangent bundle  and $Q=TM/L$ is the normal bundle
of $\mathcal F$. The metric $g_M$ determines an orthogonal
decomposition $TM=L\oplus L^\perp$, identifying $Q$ with $L^\perp$
and inducing a metric $g_Q$ on $Q$. 
Let $\nabla $ be the transverse Levi-Civita connection on $Q$,
which is torsion-free and metric with respect to $g_Q$ [\ref{Kamber}]. Let
$R^\nabla, K^\nabla,\rho^\nabla$ and $\sigma^\nabla$ be the
transversal curvature tensor, transversal sectional curvature,
transversal Ricci operator and transversal scalar curvature with
respect to $\nabla$, respectively. Let $\Omega_B^*(\mathcal F)$ be
the space of all {\it basic forms} on $M$, i.e.,
\begin{align}\label{eq1-4}
\Omega_B^*(\mathcal F)=\{\phi\in\Omega^*(M)\ | \ i(X)\phi=0,\
i(X)d\phi=0,
                     \quad \forall  X\in \Gamma L\}.
\end{align}
Then $L^2\Omega^*(M)$ is decomposed as [\ref{Lop}]
\begin{align}\label{eq1-5}
L^2\Omega(M)=L^2\Omega_B(\mathcal F) \oplus
L^2\Omega_B(\mathcal F)^\perp.
\end{align}
Now we define
the connection $\nabla $ on $\Omega _{B}^{\ast }(\mathcal F)$,
which is induced from the connection $\nabla $ on $Q $ and
Riemannian connection $\nabla ^{M}$ of $g_{M}$. This
connection $\nabla $ extends the partial Bott connection $\overset{\circ }{%
\nabla }$ given by $\overset{\circ }{\nabla }_{X}\phi =\theta(X)\phi $
for any $X\in \Gamma L$ [\ref{Kamber1}], where $\theta(X)$ is the transversal Lie derivative. Then the basic forms are
characterized by $\Omega _{B}^{\ast }(\mathcal F)={\rm
Ker}\overset{\circ }{\nabla }\subset \Gamma (\wedge Q^{\ast
}(\mathcal F))$. By a direct calculation, we have the following lemma. 
\begin{lemma} Let $(M,g_M,\mathcal F)$ be a Riemannian manifold with a foliation $\mathcal F$ and a bundle-like metric $g_M$. Then for any $X,Y, Z\in\Gamma Q$,
\begin{align*}
  [R^\nabla(X,Y),i(Z)] = i(R^\nabla(X,Y)Z).
\end{align*}
\end{lemma} 
  The exterior differential $d$ on the  de Rham complex
$\Omega^*(M)$ restricts a differential $d_B:\Omega_B^r(\mathcal
F)\to \Omega_B^{r+1}(\mathcal F)$.  
    Let $\kappa\in Q^*$ be the
mean curvature form of $\mathcal F$. Then it is well known
 that the basic part $\kappa_B$ of $\kappa$ is closed [\ref{Lop}].
 We now recall the star operator $\bar
*:\Omega^r(M)\to \Omega^{q-r}(M)$ given by
[\ref{Park},\ref{Tond1}]
\begin{align}\label{eq1-6}
 \bar * \phi =(-1)^{p(q-r)}*(\phi\wedge \chi_{\mathcal
F}),\quad\forall \phi\in\Omega^r(M),
\end{align}
where $\chi_{\mathcal F}$ is the characteristic form of $\mathcal
F$ and $*$ is the Hodge star operator associated to $g_M$. The
operator $\bar *$ maps basic forms to basic forms. 
For any $\phi,\psi\in\Omega_B^r(\mathcal F)$,  $\phi\wedge\bar
*\psi =\psi\wedge\bar *\phi$ and also ${\bar
*}^2\phi=(-1)^{r(q-r)}\phi$  [\ref{Park}]. Let $\nu$ be the transversal volume
form, i.e., $*\nu=\chi_\mathcal F$.  The pointwise inner product
$\langle\ , \ \rangle$ on $\Lambda^r Q^*$ is defined uniquely by
\begin{align}\label{eq1-8}
\langle\phi,\psi\rangle\nu=\phi\wedge \bar *\psi.
\end{align}
The global inner product $( \cdot,\cdot)_B$ on
$L^2\Omega_B^r(\mathcal F)$ is defined by
\begin{align}\label{eq1-9}
(\phi,\psi)_B = \int_M \langle\phi,\psi\rangle\mu_M,\quad
\forall\phi,\psi\in\Omega_B^r(\mathcal F),
\end{align}
where $\mu_M=\nu\wedge\chi_{\mathcal F}$ is the volume form with
respect to $g_M$.
 With
respect to this scalar product, the formal adjoint
$\delta_B:\Omega_B^r(\mathcal F)\to \Omega_B^{r-1}(\mathcal F)$ of
$d_B$ is given by [\ref{Park}]
\begin{align}\label{eq1-10}
 \delta_B\phi=(-1)^{q(r+1)+1}\bar *d_T\bar *
 \phi=\delta_T\phi+i(\kappa_B^\sharp)\phi,
 \end{align}
 where $d_T=d-\kappa_B\wedge$ and $\delta_T=(-1)^{q(r+1)+1}\bar
 *d\bar *$ is the formal adjoint operator  of $d_T$. Here $(\cdot)^\sharp$ is a $g_Q$-dual vector to $(\cdot)$. 
 The basic
Laplacian $\Delta_B$ is given by
 $ \Delta_B = d_B\delta_B
+ \delta_B d_B$.
 Let $\{E_a\} (a=1,\cdots,q)$ be a local orthonormal basic frame on $Q$.
We define $ \nabla_{\rm tr}^*\nabla_{\rm tr} :\Omega_B^r(\mathcal F)\to
\Omega_B^r(\mathcal F)$ by
\begin{align}\label{eq1-12}
\nabla_{\rm tr}^*\nabla_{\rm tr}\phi =-\sum_a \nabla^2_{E_a,E_a}
\phi+\nabla_{\kappa_B^\sharp}\phi,\quad\phi\in\Omega_B^r(\mathcal F),
\end{align}
where $\nabla^2_{X,Y}=\nabla_X\nabla_Y -\nabla_{\nabla^M_XY}$ for
any $X,Y\in TM$. Then the operator $\nabla_{\rm tr}^*\nabla_{\rm tr}$
is positive definite and formally self adjoint on the space of
basic forms [\ref{Jung}]. We define the bundle map $A_Y:\Lambda^r
Q^*\to\Lambda^r Q^*$ for any $Y\in TM$ [\ref{Kamber2}]
by
\begin{align}\label{eq1-13}
A_Y\phi =\theta(Y)\phi-\nabla_Y\phi.
\end{align}
For any $X\in\Gamma L$, $\theta(X)\phi=\nabla_X\phi$  [\ref{Kamber1}] and so $A_X\phi=0$.
Now we define the curvature endomorphism $F:\Omega_B^r (\mathcal F)\to \Omega_B^r (\mathcal F)$ by 
\begin{align}
F(\phi)=\sum_{a,b}\theta^a \wedge i(E_b)R^\nabla(E_b,
 E_a)\phi,
 \end{align}
  where  $\theta^a$ is a $g_Q$-dual 1-form to $E_a$.
Then we have
the generalized Weitzenb\"ock formula.
\begin{thm} $[\ref{Jung1}]$ On a Riemannian foliation $\mathcal F$, we have that for any $\phi\in\Omega_B^r(\mathcal
  F)$,
\begin{align*}
  \Delta_B \phi = \nabla_{\rm tr}^*\nabla_{\rm tr}\phi +
  F(\phi)+A_{\kappa_B^\sharp}\phi.
\end{align*}
 In particular, if $\phi$ is a basic 1-form, then $F(\phi)^\sharp
 =\rho^\nabla(\phi^\sharp)$.
\end{thm}
\begin{coro} On a Riemannian foliation $\mathcal F$, we have that for any $\phi\in\Omega_B^r(\mathcal F)$,
\begin{align*}
\frac12\Delta_B|\phi|^2 = \langle\Delta_B\phi,\phi\rangle -|\nabla_{\rm tr}\phi|^2 -\langle F(\phi),\phi\rangle -\langle A_{\kappa_B^\sharp}\phi,\phi\rangle.
\end{align*}
\end{coro}
Now, we recall the following generalized maximum principle.
\begin{thm} $[\ref{JLR}]$ Let  $\mathcal F$ be a Riemannian foliation on a closed, connected Riemannian manifold $(M,g_M)$. If $(\Delta_B -\kappa_B^\sharp)f\geq 0$ $($or $\leq 0)$ for any basic function $f$, then $f$ is constant.
\end{thm}
\section{The transverse conformal Killing forms}
Let $(M,g_M,\mathcal F)$ be a Riemannian manifold with a foliation $\mathcal F$ of codimension $q$ and a bundle-like metric $g_M$.
\begin{defn} {\rm A basic $r$-form $\phi\in \Omega_B^r(\mathcal F)$ is called a}  transverse conformal
Killing $r$-form {\rm if for any vector field $X\in\Gamma Q$,
\begin{align*}
\nabla_X\phi ={1\over r+1}i(X)d_B\phi -{1\over r^*+1}X^{\frak b}\wedge
\delta_T\phi,
\end{align*}
where $r^* = q-r$ and  $X^{\frak b}$ is the $g_Q$-dual 1-form of $X$. In addition, if the
basic $r$-form $\phi$ satisfies $\delta_T\phi=0$, it is called a}
transverse Killing $r$-form.
\end{defn}
Note that a transverse conformal Killing 1-form (resp. transverse
Killing 1-form) is a $g_Q$-dual form of a transversal conformal
Killing  field (resp. transversal Killing  field).

 \begin{prop}$[\ref{JK}]$ Let $\phi$ be a transverse conformal Killing $r$-form. Then
 \begin{align}
 & F(\phi) = {r\over r+1}\delta_Td_B\phi+{r^*\over r^*+1}d_B\delta_T\phi,\\
 &\nabla_{\rm tr}^*\nabla_{\rm tr}\phi={1\over r+1}\delta_Bd_B\phi +{1\over r^*+1}d_T\delta_T\phi.
 \end{align}
 \end{prop}
 \begin{lemma} $[\ref{JK}]$ Let $\phi$ be a transverse conformal
Killing $r$-form. Then
\begin{align*}
\nabla_X\nabla_Y\phi&={1\over r+1}\{i(\nabla_XY)d_B\phi +i(Y)\nabla_X d_B\phi\} \notag\\
&-{1\over r^*+1}\{\nabla_XY^{\frak b}\wedge \delta_T\phi +Y^{\frak b}\wedge
\nabla_X\delta_T\phi\}
\end{align*} 
 for any $X,Y\in\Gamma Q$. 
 \end{lemma}
   We define the operators
$R^\nabla_{\pm}(X): \wedge^r Q^* \to \wedge^{r\pm 1}Q^*$  for any
$X\in TM$ by
\begin{align}\label{eq5-1}
R^\nabla_+(X)\phi&=\sum_a \theta^a\wedge R^\nabla(X,E_a)\phi,\\
R^\nabla_-(X)\phi&=\sum_a i(E_a)R^\nabla(X,E_a)\phi.\label{eq5-2}
\end{align}
Then we have the following lemma.
\begin{lemma}\label{prop5-3} Let $\phi$ be a transverse conformal Killing
$r$-form. Then for all $X\in \Gamma Q$,
\begin{align}\label{eq5-6}
\nabla_X d_B\phi&={r+1\over r}\{R^\nabla_+(X)\phi+{1\over 
r^*+1}X^{\frak b}\wedge d_B\delta_T\phi\},\\
\nabla_X\delta_T\phi&=-{r^*+1\over r^*}\{R^\nabla_-(X)\phi
+{1\over r+1}i(X)\delta_T d_B\phi\}.\label{eq5-7}
\end{align}
\end{lemma}
{\bf Proof.} Fix $x\in M$ and  choose an orthonormal basic frame $\{E_a\}$ such that
$(\nabla E_a)_x=0$. Since $\sum_a \theta^a\wedge i(E_a)\phi = r\phi$ for any  $\phi\in\Omega_B^r(\mathcal F)$, from Lemma 3.3
\begin{align*}
R^\nabla_+(X)\phi&={r\over r+1}\nabla_X d_B\phi -{1\over r^*+1}X^{\frak b}\wedge d_B\delta_T\phi,
\end{align*}
which proves (3.5). The proof of (3.6) is similar. 
$\Box$
\begin{prop} \label{prop5-5}  Let $\phi$ be a transverse conformal Killing $r$-form. Then for any $X, Y\in \Gamma Q$,
\begin{align*}
&R^\nabla(X,Y)\phi\\
&={1\over r r^*}\Big(Y^{\frak b}\wedge i(X)-X^{\frak b}\wedge
i(Y)\Big)F(\phi)\\
&+{1\over
r}\Big(i(Y)R^\nabla_+(X)-i(X)R^\nabla_+(Y)\Big)\phi+{1\over
r^*}\Big(Y^{\frak b}\wedge R^\nabla_-(X)-X^{\frak b}\wedge R^\nabla_-(Y)\Big)\phi.
\end{align*}
\end{prop}
{\bf Proof.} Let $\phi$ be the transverse conformal Killing
$r$-form. From  Lemma 3.3,
\begin{align*}
R^\nabla(X,Y)\phi&={1\over r+1}\{i(Y)\nabla_X d_B\phi-i(X)\nabla_Y
d_B\phi\}\\
&-{1\over r^*+1}\{Y^{\frak b}\wedge \nabla_X \delta_T\phi-X^{\frak b}\wedge
\nabla_Y\delta_T\phi\}.
\end{align*}
From Lemma 3.4,
we have
\begin{align*}
&R^\nabla(X,Y)\phi\\
&={1\over r}\{i(Y)R^\nabla_+(X)-i(X)R^\nabla_+(Y)\}\phi+{1\over
r^*}\{Y^{\frak b}\wedge R^\nabla_-(X)-X^{\frak b}\wedge
R^\nabla_-(Y)\}\phi\\
&-\Big(X^{\frak b}\wedge i(Y)-Y^{\frak b}\wedge i(X)\Big)\{{1\over
r(r^*+1)}d_B\delta_T\phi+{1\over r^*(r+1)}\delta_Td_B\phi\}.
\end{align*}
Hence the proof follows  from (3.1).  $\Box$
\begin{lemma} Let $\phi$ be a transverse conformal Killing $r$-form. Then
\begin{align*}
\sum_a i(E_a)R^\nabla_-(E_a)\phi =\sum_a \theta^a\wedge R^\nabla_+(E_a)\phi=0.
\end{align*}
\end{lemma}
{\bf Proof.} Since $\phi$ is a transverse  conformal Killing $r$-form, from Proposition 3.5,
\begin{align*}
\sum_a i(E_a)R^\nabla_-(E_a)\phi&={2\over r^*}\sum_{a,b} i(E_a)i(E_b)\{\theta^b\wedge R^\nabla_-(E_a)\phi\}\\
&={2(r^*+1)\over r^*}\sum_a i(E_a)R^\nabla_-(E_a)\phi,
\end{align*}
which means that $\sum_a i(E_a)R^\nabla_-(E_a)\phi=0$. Similarly, we have
\begin{align*}
\sum_a \theta^a\wedge R^\nabla_+(E_a)\phi&={2\over r}\sum_{a,b}\theta^a\wedge\theta^b\wedge\{i(E_b)R^\nabla_+(E_a)\phi\}\\
&={2(r+1)\over r}\sum_a \theta^a\wedge R^\nabla_+(E_a)\phi,
\end{align*}
which proves the second equality. $\Box$
 
\section{ Curvatures on a K\"ahler foliation}
Let $(M,g_M,J,\mathcal F)$ be a compact Riemannian manifold with a
K\"ahler foliation $\mathcal F$ of codimension $q=2m$ and a
bundle-like metric $g_M$ [\ref{NT}]. Namely, there is a holonomy invariant almost complex structure $J:Q\to Q$ with respect to which $g_Q$ is Hermitian, i.e., $g_Q(JX,JY)=g_Q(X,Y)$ for $X,Y\in Q$ and $\nabla J=0$. Note that for any $X,
Y\in\Gamma Q$,
\begin{equation}
\Omega(X,Y)=g_Q(X,JY)
\end{equation}
defines a basic 2-form $\Omega$, which is closed as consequence of
$\nabla g_Q=0$ and $\nabla J=0$.  Then
\begin{align}
\Omega=-\frac12\sum_{a=1}^{2m}\theta^a\wedge J\theta^a.
\end{align}
  Moreover, we have the following
identities: for any $X,Y \in \Gamma Q$,
\begin{align}
 &R^{\nabla}(X,Y)J = JR^{\nabla}(X,Y),\quad
 R^{\nabla}(JX,JY) =
R^{\nabla}(X,Y).
\end{align}
Trivially, we have the following lemma.
\begin{lemma}\label{lem5-8} On a K\"ahler foliation $(\mathcal F,J)$, the following holds: 
\begin{align*}
\sum_a \theta^a\wedge \rho^\nabla(E_a)^{\frak b}=0.
\end{align*}
\end{lemma}
{\bf Proof.} By a direct calculation, we have
\begin{align*}
\sum_a\theta^a\wedge \rho^\nabla(E_a)^{\frak b}&=\sum_{a,b}\theta^a\wedge R^\nabla(E_a,JE_b)J\theta^b\\
&=\sum_{a,b,c} \theta^a\wedge g_Q(R^\nabla(E_a,JE_b)JE_b,E_c)\theta^c\\
&=\sum_{a,b}R^\nabla(E_b,JE_a)J\theta^b\wedge\theta^a\\
&=\sum_a \rho^\nabla(E_a)^{\frak b} \wedge \theta^a,
\end{align*}
which completes the proof.  $\Box$

\begin{lemma} On a K\"ahler foliation $(\mathcal F,J)$, we have that for any $\phi\in\Omega_B^r(\mathcal F)$,
\begin{align}
&\sum_a i(E_a)R^\nabla_+(E_a)\phi=\sum_a \theta^a \wedge R^\nabla_- (E_a)\phi =-F(\phi),\\
& \sum_a i(E_a)R^\nabla_-(JE_a)\phi=\sum_a\theta^a\wedge R^\nabla_+(JE_a)\phi=0.
 \end{align}
 \end{lemma}
 {\bf Proof.} The proof of (4.4) is trivial. Note that for any $X, Y \in \Gamma Q$,
\begin{align}
 R^\nabla(JX,Y) = R^\nabla(JY,X).
 \end{align}
 From (4.6), the proof of (4.5) is trivial. $\Box$
\begin{lemma} On a K\"ahler foliation $(\mathcal F,J)$, we have that for any $\phi\in\Omega_B^r(\mathcal F)$,
\begin{align*}
\sum_{a} R_+^\nabla(JE_a)i(E_a)\phi=0.
 \end{align*}
 \end{lemma}
 {\bf Proof.} Let $\phi ={1\over r!}\sum_{i_1,\cdots,i_r}\phi_{i_1\cdots i_r}\theta^{i_1}\wedge\cdots\wedge\theta^{i_r}$ be a basic $r$-form. Then by a long calculation, we have
\begin{align*}
&\sum_{a,b}\theta^a\wedge R^\nabla(JE_a,E_b)i(E_b)\phi\\
&={1\over r!}\sum_{i_1,\cdots,i_r}\sum_{a,k<l}(-1)^{k+l-1}\phi_{i_1\cdots i_r} \theta^a\wedge\{R^\nabla(JE_a,E_{i_k})\theta^{i_l}-R^\nabla(JE_a,E_{i_l})\theta^{i_k}\}\wedge\psi_{k,l}\\
&={2\over r!}\sum_{i_1,\cdots,i_r}\sum_{a,k<l}(-1)^{k+l-1} \phi_{i_1\cdots i_r}\theta^a\wedge R^\nabla(JE_a,E_{i_k})\theta^{i_l}\wedge\psi_{k,l},
\end{align*} 
where $\psi_{k,l}=\theta^{i_1}\wedge\cdots\wedge \hat \theta^{i_k}\wedge\cdots\wedge \hat\theta^{i_l}\wedge\cdots\wedge\theta^{i_r}$. 
From (4.6),
  \begin{align*}
\sum _{i_k,i_l}\phi_{i_1\cdots\i_r} R^\nabla(JE_{i_k},E_{i_l}) =0.
\end{align*}  
Hence, by the first Bianchi identity, we have
\begin{align*}
\sum_{a,i_k,i_l} \phi_{i_1\cdots i_r}\theta^a\wedge R^\nabla(JE_a,E_{i_k})\theta^{i_l}
&=\sum_{a,b,i_k,i_l}\phi_{i_1\cdots i_r} g_Q(R^\nabla(JE_a,E_{i_k})E_{i_l},E_b)\theta^a\wedge\theta^b\\
&=\sum_{a,i_k,i_l} \phi_{i_1\cdots i_r} R^\nabla(E_{i_l},E_a)J\theta^{i_k}\wedge\theta^a\\
&=\sum_{a,i_k,i_l} \phi_{i_1\cdots i_r} R^\nabla(JE_a, E_{i_k})\theta^{i_l}\wedge \theta^a \\
&=\sum_{a,i_k,i_l} \phi_{i_1\cdots i_r} R^\nabla(JE_a, E_{i_k})\theta^{i_l}\wedge\theta^a,
\end{align*}
which means 
\begin{align*}
\sum_{a,i_1,\cdots,i_r} \phi_{i_1\cdots i_r}\theta^a\wedge R^\nabla(JE_a,E_{i_k})\theta^{i_l}=0.
\end{align*}
Hence the proof is completed. $\Box$

Let $L:\Omega_B^r(\mathcal F)\to \Omega_B^{r+2}(\mathcal F)$ and $\Lambda:\Omega_b^r(\mathcal F)\to \Omega_B^{r-2}(\mathcal F)$ be given respectively by $[\ref{JJ3}]$
\begin{align}
L(\phi)=\epsilon(\Omega)\phi,\quad\Lambda(\phi)=i(\Omega)\phi,
\end{align}
where $\epsilon(\Omega)\phi=\Omega\wedge\phi$ and $i(\Omega)=-\frac12 \sum_{a=1}^{2m} i(JE_a)i(E_a)$. Trivially, for any basic forms $\phi\in\Omega_B^r(\mathcal F)$ and $\psi\in\Omega_B^{r+2}(\mathcal F)$, $\langle L(\phi),\psi\rangle =\langle\phi,\Lambda(\psi)\rangle$. Moreover, for any basic $r$-form $\phi$, $[\Lambda,L]\phi = \frac12(q-2r)\phi$.
Also, we have the following lemma.
\begin{lemma} $[\ref{JJ3}]$ On a K\"ahler foliation $(\mathcal F,J)$, we have that for any $X\in Q$,
\begin{align*}
[L,i(X)]=\epsilon(JX^{\frak b}),\ [L,\epsilon(X^{\frak b})]=[\Lambda,i(X)]=0,\ [\Lambda,\epsilon(X^{\frak b})]=-i(JX).
\end{align*}
\end{lemma}
Now, we define the operators
$\tilde J:\Omega_B^r(\mathcal F)\to \Omega_B^r(\mathcal F)$  and $S:\Omega_B^r(\mathcal F)\to \Omega_B^r (\mathcal F)$  respectively by
\begin{align} 
&\tilde J(\phi) =\sum_{a=1}^{2m}J\theta^a\wedge i(E_a)\phi,\\
&S(\phi)=\sum_{a=1}^{2m} J\theta^a \wedge i(\rho^\nabla(E_a))\phi.
\end{align}
Trivially, if $\phi\in\Omega_B^1(\mathcal F)$, then $\tilde J\phi =J\phi$. From now on, if we have no confusion, we write $\tilde J\equiv J$. 
\begin{lemma} On a K\"ahler foliation $(\mathcal F,J)$, we have that for  any $X, Y\in Q$,  
\begin{align*}
 &[J,i(X)]=i(JX),\ [J,\epsilon(X^{\frak b})]=\epsilon(JX^{\frak b}),\
 [R^\nabla(X,Y), J]=0.
 \end{align*}
 \end{lemma}
 {\bf Proof.} The first two equations are trivial. Since  $\sum_aR^\nabla(X,Y)J\theta^a\wedge i(E_a) + J\theta^a\wedge i(R^\nabla(X,Y)E_a)=0$, for any $X,Y\in Q$,
 \begin{align*}
 R^\nabla(X,Y)J\phi&= \sum_a J\theta^a \wedge i(E_a)R^\nabla(X,Y)\phi\\
 &=JR^\nabla(X,Y)\phi,
 \end{align*} 
 which proves the third equation. $\Box$ 
 \begin{lemma} On a K\"ahler foliation $(\mathcal F,J)$, we have that for any $\phi\in\Omega_B^r(\mathcal F)$,
\begin{align}
&\sum_a R^\nabla(E_a,JE_a)\phi = -2 S (\phi),\\
&\sum_{a}\theta^a\wedge R^\nabla_-(JE_a)\phi =\sum_a i(E_a)R^\nabla_+(JE_a)\phi = S(\phi).
\end{align}
\end{lemma}
{\bf Proof.}  Note that for any  $X\in \Gamma Q$,
\begin{align}
\sum_a R^\nabla (E_a,JE_a)X^{\frak b} =-2 \rho^\nabla(JX)^{\frak b}.
\end{align}
Let $\phi ={1\over r!}\sum_{i_1,\cdots,i_r}\phi_{i_1\cdots i_r}\theta^{i_1}\wedge\cdots\wedge\theta^{i_r}$. From (4.12), we have
\begin{align*}
\sum_a R^\nabla(E_a,JE_a)\phi&= -{2\over r!}\sum_{k,i_1,\cdots,i_r} \phi_{i_1\cdots i_r} \theta^{i_1}\wedge\cdots\wedge \rho^\nabla(JE_{i_k})^{\frak b}\wedge\cdots\wedge \theta^{i_r}\\
&=2\sum\theta^a\wedge i(\rho^\nabla(JE_a))\phi=-2S(\phi),
\end{align*}
which proves (4.10). From Lemma 2.1, we have
\begin{align}
\sum_a R_+^\nabla(JE_a)i(E_a)\phi&= \sum_a \theta^a \wedge R_-^\nabla(JE_a)\phi +\sum_a \theta^a \wedge i(\rho^\nabla(JE_a))\phi.
\end{align}
From Lemma 4.3 and (4.13), we have
\begin{align*}
\sum_a \theta^a\wedge R^\nabla_-(JE_a)\phi = S(\phi).
\end{align*}
Moreover, since $\sum_a R^\nabla(JE_a,E_a)\phi=\sum_a \theta^a\wedge R^\nabla_-(JE_a)\phi + \sum_a i(E_a)R^\nabla_+(JE_a)\phi$, the proof of (4.11) follows.     $\Box$
\begin{lemma} On a K\"ahler foliation $(\mathcal F,J)$, we have that for any $\phi\in\Omega_B^r(\mathcal F)$,
\begin{align}
\sum_a \theta^a\wedge J R^\nabla_-(E_a)\phi =\sum_a i(E_a)JR^\nabla_+(E_a)\phi = S(\phi) - F(J\phi).
\end{align}
\end{lemma}
{\bf Proof.}
From Lemma 4.2, Lemma 4.5 and Lemma 4.6, we have
\begin{align*}
\sum_a \theta^a\wedge JR^\nabla_-(E_a)\phi&= \sum_a J\{\theta^a\wedge R^\nabla_-(E_a)\phi\} -\sum_a J\theta^a\wedge R^\nabla_-(E_a)\phi\\
&=S(\phi)-JF(\phi) =S(\phi)-F(J\phi).
\end{align*}
The last equality in the above follows from  $[J,F]=0$.  On the other hand, from Lemma 4.5 and Lemma 4.6, we have
\begin{align*}
\sum_a i(E_a) JR^\nabla_+ (E_a)\phi &=\sum_a i(E_a) J\{\theta^b\wedge R^\nabla(E_a,E_b)\phi\}\\
&=\sum_{a,b} i(E_a)\{\theta^b\wedge JR^\nabla(E_a,E_b)\phi+J\theta^b\wedge R^\nabla(E_a,E_b)\phi\}\\
&= \sum_{a,b} i(E_a)\{\theta^b\wedge R^\nabla(E_a,E_b)J\phi+J\theta^b\wedge R^\nabla(E_a,E_b)\phi\}\\
&=-F(J\phi) +\sum_aR^\nabla(JE_a,E_a)\phi + \sum_{a} J\theta^a\wedge R^\nabla_- (E_a)\phi\\
&= S(\phi) -F(J\phi).   \quad \Box
\end{align*}

 \begin{lemma} On a K\"ahler foliation $(\mathcal F,J)$, we have
\begin{gather*}
[J,L]=[J,\Lambda]=[F,J]=[F,\Lambda]=[S,J]=[S,\Lambda]=[S,L]=0.
\end{gather*}
\end{lemma}
{\bf Proof.} From Lemma 4.5, we have
\begin{align*}
[F,J]\phi&= -\sum_{a,b}J\theta^b\wedge i(E_a) R^\nabla(E_a,E_b)\phi -\sum_{a,b}\theta^b\wedge i(JE_a) R^\nabla(E_a,E_b)\phi\\
&=0.
\end{align*}
Others are easily proved. $\Box$
 
Now, we recall the operators $d_B^c:\Omega_B^r(\mathcal F)\to\Omega_B^{r+1}(\mathcal F)$ and $\delta_B^c:\Omega_B^r(\mathcal F)\to\Omega_B^{r-1}(\mathcal F)$, which are given by $[\ref{JJ3}]$
\begin{align}
&d_B^c\phi=\sum_{a=1}^{2m}J\theta^a\wedge\nabla_{E_a}\phi,\\ 
&\delta_B^c \phi= -\sum_{a=1}^{2m}i(JE_a)\nabla_{E_a}\phi + i(J\kappa_B^\sharp)\phi.
\end{align}
 Trivially, $\delta_B^c$ is a formal adjoint of $d_B^c$ and ${\delta_B^c}^2 = {d_B^c}^2 =0$ [\ref{JJ3}]. Also, we define two operators $d_T^c$ and $\delta_B^c$ by
 \begin{align}
 d_T^c = d_B^c -\epsilon(J\kappa_B),\quad \delta_T^c = \delta_B^c -i(J\kappa_B^\sharp).
 \end{align}
Then we have the following lemma.
\begin{lemma} $[\ref{JJ3}]$ On a K\"ahler foliation $(\mathcal F,J)$, we have that 
\begin{align}
&[L,d_B]=[L,d_B^c]=0,\ \ [L,\delta_B]=-d_T^c,\ \ [L,\delta_B^c]=d_T,\\    
&[\Lambda,\delta_B]=[\Lambda,\delta_B^c]=0,\ \ [\Lambda,d_B]=\delta_T^c,\ \ [\Lambda,d_B^c]=-\delta_T,\\
&[J,d_B]=d_B^c,\ [J,\delta_B]=\delta_B^c,\
[J,d_B^c]=-d_B,\ 
 [J,\delta_B^c]=-\delta_B.
\end{align}
\end{lemma}
{\bf Proof.} Note that on K\"ahler foliations, $\nabla J=0$ and then $\nabla \tilde J=0$. Hence by Lemma 4.5, the proof follows. $\Box$
\begin{prop} On a K\"ahler foliation $(\mathcal F,J)$, we have
\begin{align}
d_T^c\delta_B + \delta_B d_T^c =d_B \delta_T^c +\delta_T^c d_B&= 0,\\
d_B^c\delta_T + \delta_T d_B^c =d_T\delta_B^c +\delta_B^c d_T&=0,\\
\delta_B\delta_B^c +\delta_B^c\delta_B =d_B d_B^c+d_B^c d_B&=0.
\end{align}
\end{prop}
{\bf Proof.} From Lemma 4.9, we have 
\begin{align*}
d_T^c\delta_B + \delta_B d_T^c = -[L,\delta_B]\delta_B -\delta_B [L,\delta_B] =0.
\end{align*}
Others are similarly proved.   $\Box$ 

Now, we put that for any $X\in TM$, 
\begin{align}
e(X)\phi=\delta_B i(X)\phi + i(X)\delta_B\phi.
\end{align}
Then we have the following.
\begin{lemma} On a K\"ahler foliation $(\mathcal F,J)$, we have that 
\begin{align*}
[J,\Delta_B]&=\theta(J\kappa_B^\sharp)+\theta(J\kappa_B^\sharp)^t,\\
[\Lambda,\Delta_B]&=e(J\kappa_B^\sharp),
\end{align*}
where $\theta(X)^t$ is a formal adjoint of $\theta(X)$ for any $X\in Q$.
\end{lemma}
Now, we recall that $\mathcal F$ is {\it minimal} if $\kappa=0$. Then we have the following corollary.
\begin{coro} On a minimal K\"ahler foliation $(\mathcal F,J)$, we have
\begin{align}
[J,\Delta_B]=[\Lambda,\Delta_B]=0.
\end{align}
\end{coro}
\section{Transverse conformal Killing forms on K\"ahler foliations}

Let $(M,g_M,J,\mathcal F)$ be a compact Riemannian manifold with a
K\"ahler foliation $\mathcal F$ of codimension $q=2m$ and a
bundle-like metric $g_M$ with respect to $\mathcal F$.
\begin{prop} On a K\"ahler foliation $(\mathcal F,J)$, if $\phi$ is a transverse conformal Killing $r$-form, then
\begin{equation}\label{5-1}
(q+r^2-qr)S(\phi)= F(J\phi).
\end{equation}
\end{prop}
{\bf Proof.} Let  $\phi$ be a transverse conformal Killing $r$-form. From Proposition 3.5,
\begin{align*}
\sum_a R^\nabla(E_a,JE_a)\phi&={2\over rr^*}J F(\phi) +{2\over r}\sum_a i(JE_a)R_+^\nabla(E_a)\phi\\
 &\ +{2\over r^*}\sum_aJ\theta^a\wedge R_-^\nabla(E_a)\phi.
\end{align*}
Hence the proof follows from Lemma 4.6. $\Box$

 From Proposition 5.1, we have the following corollary.
 \begin{coro} On a K\"ahler foliation $(\mathcal F,J)$ of codimension $q=4$, if $\phi$ is a transverse conformal Killing $2$-form, then
 \begin{align*}
 F(J\phi)=0.
 \end{align*}
 \end{coro}
 \begin{lemma} Let $\phi$ be a transverse conformal Killing $r$-form on a K\"ahler foliation. Then
 \begin{align}
 (rr^*-r-2)d_B^c\phi&=(r^*+1)d_BJ\phi -2(r+1)\delta_T L\phi,\label{5-2}\\
 (rr^*-r^*-2)\delta_T^c\phi&=(r+1)\delta_TJ\phi +2(r^*+1)d_B\Lambda\phi.\label{5-3}
 \end{align}
 \end{lemma}
 {\bf Proof.}  Since $\phi$ is a transverse conformal Killing $r$-form, from (4.7), (4.8) and (4.15), we have
 \begin{align*}
 d_B^c\phi &=\sum_a J\theta^a \wedge\nabla_{E_a}\phi={1\over r+1}Jd_B\phi -{2\over r^* +1}L\delta_T\phi.
 \end{align*}
 From the second equation in (4.18), it is trivial that $[L,\delta_T]=-d_B^c$. Hence from Lemma 4.9, we obtain
 \begin{align*}
 {rr^*-r-2\over (r+1)(r^*+1)}d_B^c\phi&={1\over r+1}d_BJ\phi -{2\over r^*+1}\delta_TL\phi,
 \end{align*} 
 which proves (\ref{5-2}). The proof of (\ref{5-3}) is similar. $\Box$

 Since $\delta_T^2\phi =-e(\kappa_B^\sharp)\phi$ for any $\phi$, from Lemma 5.3, we have that for any transverse conformal Killing $r$-form $\phi$, 
  \begin{align}
(rr^*-r-2)\delta_Td_B^c\phi&=(r^*+1)\delta_Td_BJ\phi +2(r+1)e(\kappa_B^\sharp) L\phi,\label{5-4}\\
 (rr^*-r^*-2)d_B\delta_T^c\phi&=(r+1)d_B\delta_TJ\phi,\label{5-5}\\
 (rr^*-r^*-2)\delta_T \delta_T^c\phi&= 2(r^*+1)\delta_T d_B \Lambda\phi -(r+1)e(\kappa_B^\sharp)J\phi.\label{5-5}
 \end{align}
Hence we have the following lemma.

\begin{thm} Let $(M,g_M,\mathcal F,J)$ be a closed, connected Riemannian manifold with a K\"ahler foliation of codimension $q=4$. Then for any transverse conformal Killing $2$-form,  $J\phi$ is parallel.
\end{thm}
{\bf Proof.} 
 Let $\phi$ be a transverse conformal Killing $2$-form.  Since $F(J\phi)=0$ by Corollary 5.2, we  have
\begin{align*}
d_B\delta_T J\phi + \delta_T d_B J\phi =0.
\end{align*}
Therefore, we have
\begin{align*}
\Delta_B J\phi = \theta(\kappa_B^\sharp)J\phi.
\end{align*}
Hence, by the generalized Weitzenb\"ock formula (Corollary 2.3), 
\begin{align}\label{5-7}
\frac12 (\Delta_B -\kappa_B^\sharp)|J\phi|^2 =-|\nabla_{\rm tr} J\phi |^2 \leq 0.
\end{align}
From the generalized maximum principle (Theorem 2.4), $|J\phi|$ is constant. Again, from (\ref{5-7}), we have
\begin{align*}
\nabla_{\rm tr} J\phi =0,
\end{align*}
which implies that $J\phi\in\Omega_B^2(\mathcal F)$ is parallel.     $\Box$

\begin{coro} $($cf. $[\ref{MS}])$ Let $(M,g_M,J)$ be a closed K\"ahler manifold of  dimension $4$. Then for any conformal Killing $2$-form $\phi$, $J\phi$ is parallel.
\end{coro}
On the other hand, for any basic $r$-form $\phi$, Lemma 4.9 implies that
\begin{align}
J\Lambda d_B\delta_B\phi&=d_B\delta_BJ\Lambda\phi + d_B^c\delta_B \Lambda\phi+d_B\delta_B^c\Lambda\phi +J\delta_T^c\delta_B\phi,\label{5-8}\\
J\Lambda\delta_B d_B\phi&= \delta_B d_B J\Lambda\phi + \delta_B d_B^c\Lambda\phi + \delta_B^c d_B \Lambda\phi + J\delta_B \delta_T^c \phi.\label{5-9}
\end{align}
Hence  we have the following lemma.
\begin{lemma} Let $\phi$ be a transverse conformal Killing $r(\ne q)$-form. Then
\begin{align}
J\Lambda d_B\delta_B\phi&=d_B\delta_B J\Lambda\phi + d_B\delta_B^c\Lambda\phi +d_B^c\delta_B\Lambda\phi
-{2(r^*+1)\over r^*(r+1)}J\Lambda\delta_B d_B\phi\label{5-10}\\
&\quad -Je(J\kappa_B^\sharp)\phi +{1\over r^*}J\delta_B i(\kappa_B^\sharp)J\phi,\notag\\
J\Lambda\delta_Bd_B\phi&=\delta_B d_BJ\Lambda\phi + \delta_B d_B^c\Lambda\phi +\delta_B^c d_B\Lambda\phi +{2(r^*+1)\over r^*(r+1)}J\Lambda\delta_B d_B\phi\label{5-11}\\
&\quad -{1\over r^*}J\delta_Bi(\kappa_B^\sharp)J\phi.\notag
\end{align}
\end{lemma}
{\bf Proof.}  Let $\phi$ be a transverse conformal Killing $r$-form. From (4.23), $\delta_T^c\delta_B\phi=-\delta_B\delta_T^c\phi - e(J\kappa_B^\sharp)\phi$. Hence  from Lemma 4.9 and Lemma 5.3, we have
\begin{align*}
(rr^*-r^*-2)J\delta_T^c\delta_B\phi=& -2(r^*+1) J\delta_B d_B\Lambda\phi +(r+1)J\delta_B i(\kappa_B^\sharp)J\phi\\
&-(rr^*-r^*-2)Je(J\kappa_B^\sharp)\phi\\
=&-2(r^*+1) J\Lambda\delta_B d_B\phi + 2(r^*+1)J\delta_B\delta_T^c\phi\\& -(rr^*-r^*-2)Je(J\kappa_B^\sharp)\phi+(r+1)J\delta_B i(\kappa_B^\sharp)J\phi.
 \end{align*}
Therefore, we have
\begin{align*}
r^*(r+1)J\delta_T^c\delta_B\phi= &-2(r^*+1) J\Lambda\delta_B d_B\phi -r^*(r+1)Je(J\kappa_B^\sharp)\phi\\&+ (r+1) J\delta_B i(\kappa_B^\sharp)J\phi.
\end{align*}
From  (\ref{5-8}), the proof of (\ref{5-10}) follows. The proof of (\ref{5-11}) is similar from (\ref{5-9}). $\Box$

\begin{lemma} Let $(\mathcal F,J)$ be a minimal K\"ahler foliation. Then for  a transverse conformal Killing $r$ $(2\leq r\leq q-2)$-form $\phi$,
\begin{align}
\delta_B^c d_B\Lambda\phi &= -{r+1\over (r-1)(r^*+1)}\{J\Lambda d_B\delta_B\phi +\delta_B d_B^c\Lambda \phi\},\label{5-12}\\
\delta_B d_B^c \Lambda\phi&= {r^*+1\over (r+1)(r^*-1)}\{J\Lambda \delta_B d_B\phi -\delta_B^c d_B\Lambda\phi\}.\label{5-13}
\end{align}
\end{lemma}
{\bf Proof.} From Lemma 4.9 and Proposition 4.10, we have
\begin{align}
\delta_B^cd_B\Lambda\phi&=-\Lambda d_B\delta_T^c\phi +\Lambda i(J\kappa_B^\sharp)d_B\phi -\delta_B^c\delta_T^c\phi,\label{5-14}\\
\delta_B d_B^c\Lambda\phi&=-\Lambda d_B^c\delta_T\phi + \Lambda i(\kappa_B^\sharp)d_B^c\phi +\delta_B\delta_T\phi.\label{5-15}
\end{align}
From (\ref{5-5}) and (\ref{5-14}), we have
\begin{align*}
(rr^*-r^*-2)\delta_B^c d_B\Lambda\phi=&-(r+1)\Lambda d_B\delta_T J\phi +(rr^*-r^*-2)\Lambda i(J\kappa_B^\sharp)d_B\phi\\
&-(rr^*-r^*-2)\delta_B^2\delta_T^c\phi\\
=&-(r+1)\{J\Lambda d_B\delta_T\phi-\Lambda d_B^c\delta_T\phi-\Lambda d_B\delta_T^c\phi\}\\
&+(rr^*-r^*-2)\Lambda i(J\kappa_B^\sharp)d_B\phi -(rr^*-r^*-2)\delta_B^c\delta_T^c\phi.
\end{align*}
By using (\ref{5-14}) and (\ref{5-15}), the above equation gives
\begin{align*}
\delta_B^c d_B\Lambda\phi=&-{r+1\over (r^*+1)(r-1)}\{J\Lambda d_B\delta_T\phi +\delta_Bd_B^c\Lambda\phi\}\\
&+{r+1\over (r^*+1)(r-1)}\{\Lambda i(J\kappa_B^\sharp)d_B\phi +\Lambda i(\kappa_B^\sharp)d_B^c\phi-\delta_B^c \delta_T^c\phi+\delta_B\delta_T\phi\}.
\end{align*} 
Since $\mathcal F$ is minimal, $\delta_T^c \delta_B^c\phi = \delta_B^c \delta_B^c =0$ and $\delta_T\delta_B\phi=0$. Hence the above equation proves (\ref{5-12}). From (\ref{5-15}), (\ref{5-13}) is similarly proved. $\Box$ 

Now, we put
\begin{align}
x&=J\Lambda (d_B\delta_B\phi),\quad y=J\Lambda (\delta_B d_B\phi),\quad\alpha=\delta_B^c d_B\Lambda\phi,\label{5-16}\\
 \beta&=\delta_Bd_B^c\Lambda\phi,\quad
a=d_B\delta_B J\Lambda\phi,\quad b=\delta_B d_B J\Lambda\phi.\label{5-17}
\end{align}
From now on, we assume that $\mathcal F$ is minimal. From Lemma 5.6, we have 
\begin{align}
x&= a-\alpha-\beta-{2(r^*+1)\over r^*(r+1)}y,\label{5-18}\\
y&= b + \alpha +\beta + {2(r^*+1)\over r^*(r+1)}y.\label{5-19}
\end{align}
Hence from (\ref{5-18}) and (\ref{5-19}), we have
\begin{align}
(rr^*-r^*-2)y&=r^*(r+1) (b +\alpha +\beta),\label{5-20}\\
(rr^*-r^*-2)x&=(rr^*-r^*-2)a -2(r^*+1)b -r^*(r+1)(\alpha+\beta).\label{5-21}
\end{align}
On the other hand, from Lemma 5.7, we have
\begin{align}
\alpha&= -{r+1\over (r^*+1)(r-1)}(x+\beta),\label{5-22} \\
\beta&= {r^*+1\over (r+1)(r^*-1)}(y-\alpha).\label{5-23}
\end{align}
Note that $rr^*-r^*-2=0$ if and only if $q=4$. Hence from (\ref{5-20}), (\ref{5-21}), (\ref{5-22}) and (\ref{5-23}), if $q\ne 4$, then 
\begin{align}
&\lambda_1\lambda_3 b = (1-\lambda_1 \lambda_3)\beta +\lambda_3 (1-\lambda_1)\alpha,\label{5-24}\\
&\lambda_2 a + \lambda_2(1-\lambda_1)b = (\lambda_1\lambda_2 -1)\alpha +\lambda_2(\lambda_1-1)\beta,\label{5-25}
\end{align}
where $\lambda_1={r^*(r+1)\over rr^*-r^*-2},\ \lambda_2={r+1\over (r^*+1)(r-1)}$ and $\lambda_3 = {r^*+1\over (r+1)(r^*-2)}$.
Hence we have the following theorem.
\begin{thm} Let $(M,g_M,J,\mathcal F)$ be a closed Riemannian manifold with a minimal K\"ahler foliation $\mathcal F$ of codimension $q=2m$ and a bundle-like metric $g_M$.  Then for any transverse conformal Killing $r\ (2\leq r\leq q-2)$-form $\phi$, $J\Lambda\phi$ is basic-harmonic.
\end{thm}
{\bf Proof.}
From  Lemma 4.9, $d_B\delta_B^c + \delta_B^c d_B = \theta(J\kappa_B^\sharp)\phi$. Hence we have
\begin{align*}
\int_M \langle b,\alpha\rangle\mu_M = \int_M \langle d_B J\Lambda\phi,\theta(J\kappa_B^\sharp)d_B\Lambda\phi\rangle\mu_M.
\end{align*}
Since $\mathcal F$ is minimal, we have
\begin{align}\label{5-26}
\int_M \langle b,\alpha\rangle\mu_M=0.
\end{align}
Similarly, we have
\begin{align}
\int_M \langle \beta,\alpha\rangle\mu_M =0\label{5-27}
\end{align}
and
\begin{align}
\int_M \langle a,b\rangle\mu_M =\int_M \langle a,\beta\rangle\mu_M =0.\label{5-28}
\end{align}
(i) In case of $q\ne 4$. From  (\ref{5-24}), (\ref{5-26}) and (\ref{5-27}), we have 
\begin{align*}
\lambda_3(1-\lambda_1)\int_M |\alpha|^2\mu_M=0.
\end{align*}
Since $\lambda_3\ne 0 $ and $\lambda_1\ne 1$, $\alpha=0$. 
From (\ref{5-20}) and (\ref{5-28}), $a=0$. Therefore, from (\ref{5-24}) and (\ref{5-25}), since $\lambda_2(1-\lambda_1)\ne 0$, we have
\begin{align*}
\lambda_1\lambda_3 b = (1-\lambda_1 \lambda_3)\beta,\quad b=-\beta.
\end{align*}
 Hence $b=\beta=0$. Therefore, $x=y=0$. So from Corollary 4.12, $\Delta_B J\Lambda\phi=J\Lambda\Delta_B\phi=x+y=0$. That is, $J\Lambda\phi$ is basic-harmonic. (ii) In case of $q=4$. From Theorem 5.4, $J\phi\in\Omega_B^2(\mathcal F)$ is parallel and so basic-harmonic, i.e., $\Delta_BJ\phi=0$. Hence from  Corollary 4.12, $\Delta_BJ\Lambda\phi=0$, i.e., $J\Lambda\phi$ is  basic-harmonic.  
    $\Box$

\begin{coro} Let $(M,g_M,J,\mathcal F)$ be as in Theorem 5.8. Then for a transverse conformal Killing $r\ (2\leq r\leq q-2)$-form $\phi$, $J\Lambda \phi$ is parallel.
\end{coro}
{\bf Proof.} Let $\phi$ be a transverse conformal Killing form. Since $\mathcal F$ is minimal, from Theorem 5.8, $\Delta_B (J\Lambda\phi)=0$. Hence from the generalized Weitzenb\"ock formula (Theorem 2.2), we have
\begin{align*}
F(J\Lambda\phi) +\nabla_{\rm tr}^*\nabla_{\rm tr} J\Lambda \phi =0.
\end{align*}
On the other hand, from  Proposition 3.2, we have
\begin{align}
F(J\Lambda\phi)  = {r\over r+1} y + {r^*\over r^*+1} x.
\end{align}
In the proof of Theorem 5.8, $x=y=0$. Hence $F(J\Lambda\phi)=0$, which means  that $J\Lambda \phi$ is parallel. $\Box$

\begin{thm}  Let $(M,g_M,J,\mathcal F)$ be a closed Riemannian manifold with a minimal  K\"ahler foliation $\mathcal F$ of codimension $q=2m(\ne 4)$ and a bundle-like metric $g_M$. Then for a transverse conformal Killing $r (r\ne m,\ 2\leq r\leq q-2)$-form,  $J\phi$ is parallel. 
\end{thm}
{\bf Proof.}
 Let $\phi$ be a transverse conformal Killing $r$-form. Then $\bar *\phi$ is also a transverse conformal Killing $(q-r)$-form [\ref{JK}]. Hence by Corollary 5.9, $J\Lambda  \bar *\phi$ is parallel. Since $[\nabla_{\rm tr},\bar *]=0$, $[J,\bar *] =0$ and $L\bar *=\bar *\Lambda$, $\bar * J\Lambda \bar * \phi = \pm LJ\phi$ is parallel.
Note that $(m-r)J\phi = [\Lambda,L]J\phi$.  Since $[L,\nabla_{\rm tr}]=[\Lambda,\nabla_{\rm tr}]=[J,\Lambda]=0$, from Corollary 5.9, we get
\begin{align*}
(m-r)\nabla_{\rm tr} J\phi =\nabla_{\rm tr}\Lambda LJ\phi -\nabla_{\rm tr}L\Lambda J\phi=\Lambda\nabla_{\rm tr}LJ\phi -L\nabla_{\rm tr}J\Lambda\phi=0.
\end{align*}
Hence  if $r\ne m$, then $J\phi$ is parallel.  $\Box$

\bigskip
\noindent
{\bf Remark.}  (1)  When $q=4$,   $J\phi$ is parallel for any transverse conformal 2-form $\phi$ (Theorem 5.4). 

(2)  For the point foliation, Theorem  5.10  has been proved in [\ref{MS}].

\bigskip
\noindent{\bf Question}. When $\mathcal F$ is not minimal, is Theorem 5.10 true? 

\bigskip
\noindent{\bf Acknowledgements.} The author would like to thank the referee for his or her significant corrections and kind comments. This research was supported by  the  Basic Science Research Program
       through the National Research Foundation of Korea(NRF) funded
       by the Ministry of Education, Science and Technology(2010-0021005).

\vskip 1.0cm Seoung Dal Jung

\noindent Department of Mathematics and Research Institute for basic Sciences, Jeju National University, Jeju 690-756, Korea(e-mail: sdjung@jejunu.ac.kr)

\end{document}